\documentclass[12pt]{article}
\usepackage{arxiv}

\usepackage{siunitx}
\usepackage[utf8]{inputenc} 
\usepackage[T1]{fontenc}    
\usepackage{amsmath}
\usepackage{bm}
\usepackage{booktabs}       
\usepackage{amsfonts}       
\usepackage{amssymb}
\usepackage{nicefrac}       
\usepackage{microtype}      
\usepackage{soul}
\usepackage[table]{xcolor}
\definecolor{medium_blue}{rgb}{0, 0, 0.804}
\usepackage{hyperref}
\hypersetup{colorlinks=true,linkcolor=medium_blue, citecolor=medium_blue}
\usepackage{leftidx}
\usepackage{graphicx}
\usepackage{threeparttable}
\usepackage{multirow}
\usepackage{float}
\usepackage{physics}
\usepackage{lscape}
\usepackage{subcaption}

\title{Cardiac Electrophysiology Meshfree Modeling through the Mixed Collocation Method}

\author{
  Konstantinos A. ~Mountris\thanks{[mail] kmountris@unizar.es \quad [url] https://www.mountris.org} \\
  Arag\'on Institute of Engineering Research, IIS Arag\'on\\
  CIBER-BBN\\
  University of Zaragoza\\
  Spain, Zaragoza, ZGZ 50018 \\
  \texttt{kmountirs@unizar.es} \\
  \And
  Esther ~Pueyo \\
  Arag\'on Institute of Engineering Research, IIS Arag\'on, \\
  CIBER-BBN\\
  University of Zaragoza\\
  Spain, Zaragoza, ZGZ 50018 \\
  \texttt{epueyo@unizar.es} \\
}

\begin{document}
\maketitle

\begin{abstract}
We present the meshfree Mixed Collocation Method (MCM) to solve the monodomain model for numerical simulation of cardiac electrophysiology. We apply MCM to simulate cardiac electrical propagation in 2D tissue sheets and 3D tissue slabs as well as in realistic large-scale biventricular anatomies. Capitalizing on the meshfree property of MCM, we introduce an immersed grid approach for automated generation of nodes in the modeled cardiac domains. We demonstrate that MCM solutions are in agreement with FEM solutions, thus confirming their suitability for simulation of cardiac electrophysiology both in healthy and disease conditions, such as left-bundle-branch block (LBBB) and myocardial infarction. Despite the fact that the computational time for MCM calculations is longer than for FEM, its efficiency in dealing with domains presenting irregularity, nonlinearity and discontinuity make MCM a promising alternative for heart's electrical investigations.
\end{abstract}

\keywords{meshfree \and mixed collocation method \and MCM \and cardiac electrophysiology \and monodomain model}

\section{Introduction}
In the last decades, computational modeling and simulation has taken a growing role as a method to deepen the understanding of cardiac function in health and disease \cite{lopez2015three,niederer2019computational}. Novel \textit{in silico} models of increasing complexity are continuously being developed to simulate the electrophysiology \cite{sampedro2020characterization} and mechanics of the heart \cite{rama2018towards} from the cell \cite{pueyo2016experimentally} to the whole-organ level \cite{chabiniok2016multiphysics}. At the tissue and organ levels, electrophysiology is simulated by using the well-known bidomain \cite{tung1978bidomain} and monodomain \cite{keener2009mathematical} models. The latter is a simplified version of the former under the assumption of equal anisotropy ratios for the intracellular and extracellular spaces. The monodomain model is more computationally efficient than its bidomain counterpart and is able to produce accurate transmembrane potential values in the absence of extracellularly applied currents \cite{potse2006comparison}.

Commonly, state-of-the-art simulators \cite{mirams2013chaste,vigmond2003carp} employ the Finite Element Method (FEM) to solve either the bidomain or the monodomain model for the simulation of cardiac electrophysiology. Despite the fact that FEM is a mature and robust numerical method, its requirement for a good-quality mesh poses challenges to generate realistic heart models with reasonable computational cost especially when cardiac mechanics are also accounted for and large deformations are aimed to be simulated \cite{lluch2019breaking}.

Alternative meshfree methods can alleviate the mesh requirement and have been proposed for both cardiac electrophysiological \cite{lluch2017smoothed,zhang2012meshfree,mountris2019novel} and mechanical \cite{lluch2019breaking,legner2014studying} simulations. Among the different proposed meshfree solutions, models based on the Element-Free Galerkin (EFG) method offer high convergence rate and high resolution of localized steep gradients \cite{belytschko1994element}. Nevertheless, special treatment for the imposition of essential boundary conditions is required since the approximation functions do not possess the Kronecker delta property. 
Recently, Cell-based Maximum Entropy (CME) approximants were used in EFG to alleviate this limitation \cite{mountris2020cell}. 

CME possesses the weak Kronecker delta property where approximation functions of internal nodes vanish on the boundaries. Therefore, essential boundary conditions can be imposed directly as in FEM. However, the CME approximants give rise to complex integrals requiring a large number of quadrature points for accurate integration that may lead to increased computational cost. Similarly, the computational cost of methods based on the Smoothed Particle Hydrodynamics (SPH) may be significantly higher than that of mesh-based methods. Furthermore, the standard SPH formulation may imply inaccurate computation of gradients of constant and linear fields (first-order incomplete approximation) \cite{chen1999corrective}. To overcome these problems, the total Lagrangian formulation of SPH \cite{bonet1999variational} and gradient normalization \cite{chen1999corrective} were applied in \cite{lluch2017smoothed,lluch2019breaking} to accurately simulate the propagation of the electrical impulse in the heart and cope with large deformations in the context of cardiac mechanics.

In the present study, we propose the Mixed Collocation Method (MCM) as an alternative to mesh-based and previously mentioned meshfree methods for cardiac electrophysiology simulation. MCM is based on the Meshfree Local Petrov-Galerkin (MLPG) method \cite{atluri2004meshless,atluri2005basis}. MLPG implies quadrature over locally-defined domains providing the flexibility to select the trial and test functions from different spaces. In the mixed collocation variant of the MLPG method, the Dirac function is used as test function and both the field function and its gradient are interpolated by the trial function. As a result, the computational cost is decreased since the local integration is reduced to nodal summation. Moreover, while collocation methods suffer from inaccurate imposition of Neumann boundary conditions \cite{libre2008}, the accuracy is ameliorated in MCM due to the interpolation of the field function's gradient. In the seminal work on MCM \cite{atluri2006}, the moving least squares (MLS) approximation \cite{lancaster1981surfaces} was used as the trial function. Recently, the radial point interpolation (RPI) \cite{wang2002point} was proposed as an alternative to MLS in MCM \cite{mountris2020radial}. It was demonstrated that accuracy is improved when RPI trial functions replace MLS. This was mainly attributed to the delta Kronecker property of RPI that allows the direct imposition of essential boundary conditions in contrast to MLS. 

Here we investigate the application of MCM with interpolating trial functions for the solution of the monodomain equation in a series of cardiac electrophysiology simulations in 2D and 3D domains. We evaluate  RPI as well as moving Kriging interpolation (MKI) \cite{gu2003moving} and we compare the obtained solutions with state-of-the-art FEM solution. The structure of the paper is the following. In section \ref{sec:mcm_theory}, we derive the form of the cardiac monodomain equation using the MCM method. In section \ref{sec:interpolation}, we describe the mathematical formulation of RPI and MKI interpolations. In section \ref{sec:numerical}, we evaluate the solution of the monodomain model with the MCM method in several 2D and 3D problems and we report on the accuracy and efficiency of the method in comparison to FEM. Finally, in section \ref{sec:conclusions}, we discuss some concluding remarks.

\section{Mixed collocation form of the monodomain equation} \label{sec:mcm_theory}

The propagation of an electrical impulse in the heart is modeled through the monodomain model by the reaction-diffusion equation:
\begin{equation} \label{eq:monodomain_pde}
\begin{array}{ll}
\partial V(t) / \partial t = \bm{\nabla} \cdot (\bm{D\nabla}V(t)) -I_{ion}(V(t))/C & \textrm{ in } \Omega \\
\bm{n} \cdot (\bm{D\nabla}V(t)) = 0 & \textrm{ on } \partial\Omega
\end{array}
\end{equation}

\noindent where $\partial V / \partial t$ is the time derivative of the transmembrane potential $V$, $I_{ion}$ is the total ionic current and $C$ is the cell capacitance per unit surface area. $\Omega$ and $\partial \Omega$ denote the domain of interest and its boundary, respectively, and $\bm{n}$ is the outward unit vector normal to the boundary. $\bm{D}$ is the diffusion tensor given by:
\begin{equation} \label{eq:diffusion_tensor}
\bm{D} = d_0 \left[(1-\rho)\bm{f}\otimes\bm{f} + \rho \bm{I} \right]
\end{equation}

\noindent where $d_0$ denotes the diffusion coefficient along the cardiac fibers, $\rho \leq 1$ is the transverse-to-longitudinal conductivity ratio, $\bm{f}$ denotes the cardiac fiber direction vector, $\bm{I}$ is the identity matrix and $\otimes$ is the tensor product operation.

The transmembrane potential $V$ is characterized by a rapid upstroke phase followed by a slow repolarization period, thus requiring to solve a system of stiff ordinary differential equations (ODEs) to reproduce it accurately \textit{in silico}. Therefore, it is common for large scale tissue simulations to decouple Equation (\ref{eq:monodomain_pde}) using the operator-splitting method \cite{qu1999advanced} to improve computational efficiency. By applying operator-splitting, the decoupled monodomain equation is given by:

\begin{equation} \label{eq:decoupled_monodomain_pde}
\begin{array}{ll}
\partial V(t) / \partial t = -I_{ion}(V(t))/C & \textrm{ in } \Omega \\
\partial V(t) / \partial t = \bm{\nabla} \cdot (\bm{D\nabla}V(t)) & \textrm{ in } \Omega \\
\bm{n} \cdot (\bm{D\nabla}V(t)) = 0 & \textrm{ on } \partial\Omega
\end{array}
\end{equation}

\subsection{Deriving the mixed collocation form} \label{subsec:mcm_formulation}

To derive the mixed collocation form of the monodomain model we consider the discretization of the domain into a cloud of $N$ arbitrarily distributed field nodes. We write the diffusion term from Equation (\ref{eq:decoupled_monodomain_pde}) for each field node $I$ in the set of $N$ nodes as:
\begin{equation} \label{eq:diffusion_term}
\partial V_I(t) / \partial t = \bm{\nabla} \cdot \bm{q}_I(t) \textrm{ in } \Omega \\
\end{equation}

\noindent where
\begin{equation} \label{eq:potential_flux}
\bm{q}_I(t) = \bm{D}_I \bm{\nabla}V_I(t)
\end{equation}

\noindent denotes the transmembrane potential flux at the field node $I$. Interpolating the transmembrane potential and the transmembrane potential flux we obtain:
\begin{equation} \label{eq:voltage_interp}
    V_I(t) = \sum_{i=1}^n \phi_I^i V_I^i(t)
\end{equation}
\begin{equation} \label{eq:gradient_interp}
    \bm{q}_{I}(t) = \sum_{i=1}^n \phi_I^i \bm{q}_{I}^i(t)
\end{equation}

\noindent where $n$ is the number of field nodes in the local support domain of the node $I$. $\phi_I^i$ is the $i^{th}$ component of the vector $\bm{\phi}_I$ containing the basis function of the meshfree approximation at each of the $n$ nodes in the local support domain of $I$. $V_I^i(t)$ denotes the transmembrane potential at node $i$ of the local support domain of $I$ and $\bm{q}_{I}^{i}(t)$ the corresponding transmembrane potential flux vector.

By introducing Equation (\ref{eq:voltage_interp}) into Equation (\ref{eq:potential_flux}), we can express the transmembrane potential flux in terms of the nodal transmembrane potential as follows:
\begin{equation} \label{eq:gradient_nodal_voltage}
\bm{q}_I(t) = \bm{D}_I \sum_{i=1}^n \bm{\nabla}\phi_I^i V_I^i(t), \quad I = 1,2, \dotso, N
\end{equation}

\noindent or in matrix form and removing the dependence on $t$ for simplicity:
\begin{equation} 
    \bm{q}_I = \bm{K}_{a,I} \bm{V}_I, \quad I = 1,2, \dotso, N.
\end{equation}

\noindent where $\bm{K}_{a,I}$ contains the spatial partial derivatives of the meshfree basis function at the nodes in the support domain of node $I$ scaled by $\bm{D}_I$. 


By grouping the equations for all field nodes $I$, with $I = 1,2, \dotso, N$, the following equation in matrix form can be obtained:
\begin{equation} \label{eq:nodal_voltage_flux_matform}
    \bm{q} = \bm{K_a} \bm{V} .
\end{equation}


Finally, introducing Equations (\ref{eq:voltage_interp}--\ref{eq:gradient_nodal_voltage}) into Equation (\ref{eq:diffusion_term}), the mixed collocation formulation of the monodomain model's diffusion term can be obtained in terms of the transmembrane potential:
\begin{equation} \label{eq:mc_diffusion_nodal_voltage}
\sum_{i=1}^n \phi_I^i \partial V_I^i(t) / \partial t - \sum_{i=1}^n \div{( \bm{D}_I\grad{\phi_I^i}V_I^i(t))} = 0, \quad I = 1,2, \dotso ,N .
\end{equation}

\noindent By grouping the equations for all nodes $I$ and removing the dependence on $t$, as above, Equation (\ref{eq:mc_diffusion_nodal_voltage}) can be written in the equivalent matrix form:

\begin{equation} \label{eq:mc_diffusion_nodal_voltage_matfom}
    \bm{M} \dot{\bm{V}} + \bm{KV} = \bm{0}, \quad \bm{K} = \bm{K_sK_a}
\end{equation}

\noindent where $\bm{M}$ is the sparse matrix collecting the basis functions and $\bm{K}$ denotes the stiffness matrix.

\subsection{Boundary conditions imposition} \label{subsec:mcm_bc}

For the monodomain model, the domain $\Omega$ is assumed to be isolated in the sense that no current can flow in or out of the boundary $\partial \Omega$. To model electrical isolation, we enforce the Neumann boundary conditions (BC) in mixed collocation using the penalty method described in \cite{overvelde2012moving}. From Equations (\ref{eq:monodomain_pde}) and (\ref{eq:potential_flux}), the Neumann BC imposition on the $\gamma_{bc}$ nodes of the Neumann boundary $\partial \Omega$ at a given time $t$, where $\gamma_{bc} \subset \{ 1, \cdots, N\}$, can be written in matrix form as:

\begin{equation} \label{eq:flux_bc_mat}
    \bm{N}_{bc} \bm{q}_{bc} = \bm{0},    
\end{equation}

\noindent where $\bm{q}_{bc}$ is the vector collecting the transmembrane potential fluxes at $\gamma_{bc}$ nodes and $\bm{N}_{bc}$ is the matrix containing the normal vectors given by:

\begin{equation}
    \bm{N}_{bc} = \begin{bmatrix}
                    \bm{n}_1 & & 0\\
                    & \ddots & \\
                    0 & & \bm{n}_{\gamma_{bc}}
                \end{bmatrix}
\end{equation}

\noindent The Neumann BCs are enforced at the $\gamma_{bc}$ nodes by multiplying Equation (\ref{eq:flux_bc_mat}) with the penalty factor $\alpha \bm{N}_{bc}^T$ and adding it to Equation (\ref{eq:nodal_voltage_flux_matform}) to obtain:

\begin{equation} \label{eq:modified_nodal_voltage_flux}
    \bm{q}_{bc} + \alpha \bm{N}_{bc}^T \bm{N}_{bc} \bm{q}_{bc} = \bm{K}_a^{bc} \bm{V}^{bc}.
\end{equation}

\noindent By rearranging terms, Equation (\ref{eq:modified_nodal_voltage_flux}) can be written as:

\begin{equation} \label{eq:neumann_imposition}
  \bm{q}_{bc} = \{\bm{I} + \alpha \bm{N}_{bc}^T \bm{N}_{bc}\}^{-1} \{\bm{K_a}^{bc} \bm{V}^{bc} \} = \bm{Q}^{-1} \{\bm{K}_a^{bc} \bm{V}^{bc} \},
\end{equation}

\noindent where $\bm{I}$ is the identity matrix and $\bm{Q} = \bm{I} + \alpha \bm{N}_{bc}^T \bm{N}_{bc}$. Combining Equations (\ref{eq:mc_diffusion_nodal_voltage_matfom}) and (\ref{eq:neumann_imposition}), the matrix form of the monodomain model's diffusion term is given by:

\begin{equation} \label{eq:heat_balance_nodal_mat_mod}
    \bm{M} \dot{\bm{V}} + \bm{K'V} = \bm{0},
\end{equation}

\noindent where $\quad \bm{K'} = \left[\bm{K_s^{bc} Q^{-1} K_a^{bc}} \quad , \quad \bm{K_s^{in} K_a^{in}}\right]$. Here, the superscripts $bc$ and $in$ connote the row entries of the matrices $K_s$ and $K_a$ that correspond to the $\gamma_{bc}$ nodes on $\partial \Omega$ and the $\gamma_{in}$ nodes in $\Omega$, respectively, such that $\gamma_{bc} \cup \gamma_{in} = \{ 1, \cdots, N\}$. The value of the penalty factor $\alpha$ should be sufficiently large to ensure accurate imposition of the boundary condition, as  instability issues may arise if $\alpha$ is too large. In this study, we used $\alpha = 10^6$ as it was found to be the optimal value in \cite{mountris2020radial}.

\section{Interpolating meshfree approximants} \label{sec:interpolation}

One of the advantages of MCM, being a meshfree method, is the flexibility that offers on the choice of the trial function $\bm{\phi}$. In this work, we consider only trial functions that possess the delta Kronecker property, namely the radial point interpolation (RPI) \cite{wang2002point} and the moving Kriging interpolation (MKI) \cite{gu2003moving}.

\subsection{Radial point interpolation}

The RPI trial function $\bm{\phi}_I$ for a field node $I$ is obtained by:

\begin{equation} \label{eq:rpi}
    \bm{\phi}_I = \{ \bm{r}_I \; \; \bm{p}_I\} \bm{G}^{-1},
\end{equation}

\noindent where $\bm{r}_I$ is the radial basis function (RBF) for node $I$. For a node $i$ in the support domain of $I$ (with $I$ included in its support domain), $\bm{r}_i$ is given by:

\begin{equation} \label{eq:rbf}
    \bm{r}_i = \left[ r_{i1} \; r_{i2} \; \dotso \; r_{in} \right].
\end{equation}

\noindent Different RBFs, such as multiquadric, Gaussian, etc., can be used. In this work, we used the multiquadric RBF (MQ-RBF). The value of the MQ-RBF for the 3D case, is calculated as:

\begin{equation} \label{eq:mq_rbf}
    r_{Ii} = \left(d^{i2}_{I} + r_c^2\right)^q = \left[ (x_I-x_i)^2 + (y_I-y_i)^2 + (z_I-z_i)^2 + r_c^2  \right]^q
\end{equation}

\noindent where $d^{i}_{I}$ denotes the Euclidean distance between nodes $i$ and $I$, $i=1,\dotso,n$ and $r_c$ and $q$ are positive-valued shape parameters of the MQ-RBF. For spherical support domains, the shape parameter $r_c$ is given by:

\begin{equation} \label{eq:rpi_param_rc}
    r_c = \alpha_c d_c
\end{equation}

\noindent where $d_c$ denotes the radius of the support domain of node $I$ and $\alpha_c$ is a dimensionless constant. RBF fail to reconstruct exactly a linear polynomial field, therefore, the RPI is enriched with the linear polynomial basis $\bm{p_I}$ to ensure $C^1$ continuity. For 3D problems, $\bm{p}_I$ is given by:

\begin{equation} \label{eq:linear_poly}
    \bm{p}_I = [1 \; x_I \; y_I \; z_I].
\end{equation}

\noindent Finally the matrix $\bm{G}$ is given by:

\begin{equation}
    \bm{G} = \left[
             \begin{matrix}
                \bm{R}   & \bm{P}\\
                \bm{P}^T & \bm{0}
                \end{matrix}
                \right]_{(n+m \times n+m)},
\end{equation}

\noindent where $m$ is the number of components of the polynomial basis ($m=4$ for linear $\bm{p}_I$ in 3D). $\bm{R}$ and $\bm{P}$ denote the RBF and polynomial basis moment matrices:

\begin{equation} \label{eq:moment_mats}
\bm{R} = \left[
\begin{matrix}
r_{11} & r_{12} & \hdots & r_{1n}\\
r_{21} & r_{22} & \hdots & r_{2n}\\
\hdots & \hdots & \hdots & \hdots\\
r_{n1} & r_{n2} & \hdots & r_{nn}\\
\end{matrix}
\right]_{(n \times n)} , \; \;  \bm{P} = \left[
                              \begin{matrix}
                              1 & x_1 & y_1 & z_1\\
                              1 & x_2 & y_2 & z_2\\
                              \hdots & \hdots & \hdots & \hdots\\
                              1 & x_n & y_n & z_n\\
                              \end{matrix}
                              \right]_{(n \times m)}.
\end{equation}

\subsection{Moving Kriging Interpolation}
The moving Kriging interpolation (MKI) has similar interpolation properties to RPI but it does not require polynomial enrichment to ensure $C^1$ continuity. The trial function $\bm{\phi}$ at node $I$ is given by:

\begin{equation} \label{eq:mki}
    \bm{\phi}_I = \bm{p}_I\bm{A} + \bm{c}_I\bm{B} 
\end{equation}

\noindent where $\bm{p}_I$ is the linear polynomial basis defined in Equation (\ref{eq:linear_poly}) and $\bm{c}_I$ denotes the correlation function for node $I$. For a node $i$ in the support domain of $I$, $\bm{c}_i$ is given by:
\begin{equation} \label{eq:correlation}
    \bm{c}_i = \left[ c_{i1} \; c_{i2} \; \dotso \; c_{in} \right]
\end{equation}

\noindent where $c_{ij}$ is the value of the correlation function $\bm{c}_i$ at the $j^{th}$ node of the support domain. In this work, we use the MQ-RBF as the correlation function (Equation \ref{eq:mq_rbf}). The matrices $\bm{A}$ and $\bm{B}$ are obtained by:

\begin{align} \label{eq:mki_mats}
\begin{split}
    \bm{A} &= (\bm{P}^T\bm{C}^{-1}\bm{P})^{-1} \bm{P}^T \bm{C}^{-1} \\
    \bm{B} &= \bm{C}^{-1}(\bm{I} - \bm{P}\bm{A})
\end{split}
\end{align}

\noindent where $\bm{I}$ is the $n \times n$ identity matrix, $\bm{P}$ is the $n \times m$ moment matrix of the linear polynomial basis given by Equation (\ref{eq:moment_mats}) and $\bm{C}$ is the $n \times n$ correlation matrix for the $n$ nodes in the support domain of $I$ given by:

\begin{equation} \label{eq:moment_mats}
\bm{C} = \left[
\begin{matrix}
c_{11} & c_{12} & \hdots & c_{1n}\\
c_{21} & c_{22} & \hdots & c_{2n}\\
\hdots & \hdots & \hdots & \hdots\\
c_{n1} & c_{n2} & \hdots & c_{nn}\\
\end{matrix}
\right]_{(n \times n)} .
\end{equation}

\section{Numerical examples} \label{sec:numerical}

In this section, we investigate the accuracy and efficiency of MCM using both RPI and MKI as trial functions. We consider regular and irregular nodal distributions in 2D tissue sheets and 3D tissue slabs as well as in realistic anatomical models and we perform a comparison of MCM and FEM simulation results. MQ-RBFs for both RPI and MKI trial functions are constructed using $\alpha_c = 1.03$ and $q=1.42$ in 2D simulations. In 3D, they are constructed using $\alpha_c = 1.03$ and $q=1.82$, as these combinations of parameters are found to minimize the difference with FEM. In all examples, time integration is performed explicitly using the dual adaptive explicit time integration (DAETI) algorithm \cite{mountris2020dual} with adaptive time step $dt = 0.1$ ms. Human ventricular cellular electrophysiology is represented by the O'Hara et al. cell model \cite{ohara2011simulation}. Simulations are performed on a laptop with Intel\textsuperscript{\textregistered} Core\texttrademark i7-4720HQ CPU and 16 GB of RAM.

\subsection{Electrical propagation in a 2D tissue sheet}

We consider a $5 \times 5$ cm human ventricular tissue sheet, where transmural heterogeneities are included by defining endocardial, midmyocardial and epicardial regions in a 5:2:3 ratio. The cardiac fiber direction vector $\bm{f}$ is considered to be parallel to the $x$-axis. We use a diffusion coefficient value $d_0 = 0.0013$ cm$^2$/ms and a transverse-to-longitudinal conductivity ratio $\rho = 0.2$. Stimuli with amplitude twice the diastolic threshold, period $t_T = 1$ s and duration $t_d = 1$ ms are applied on the left side of the tissue ($x = 0$ cm). The propagation of the action potential (AP) is simulated for a total time $t_s = 3$ s. 

\begin{figure}[H]
    \centering
    \includegraphics[width=\textwidth]{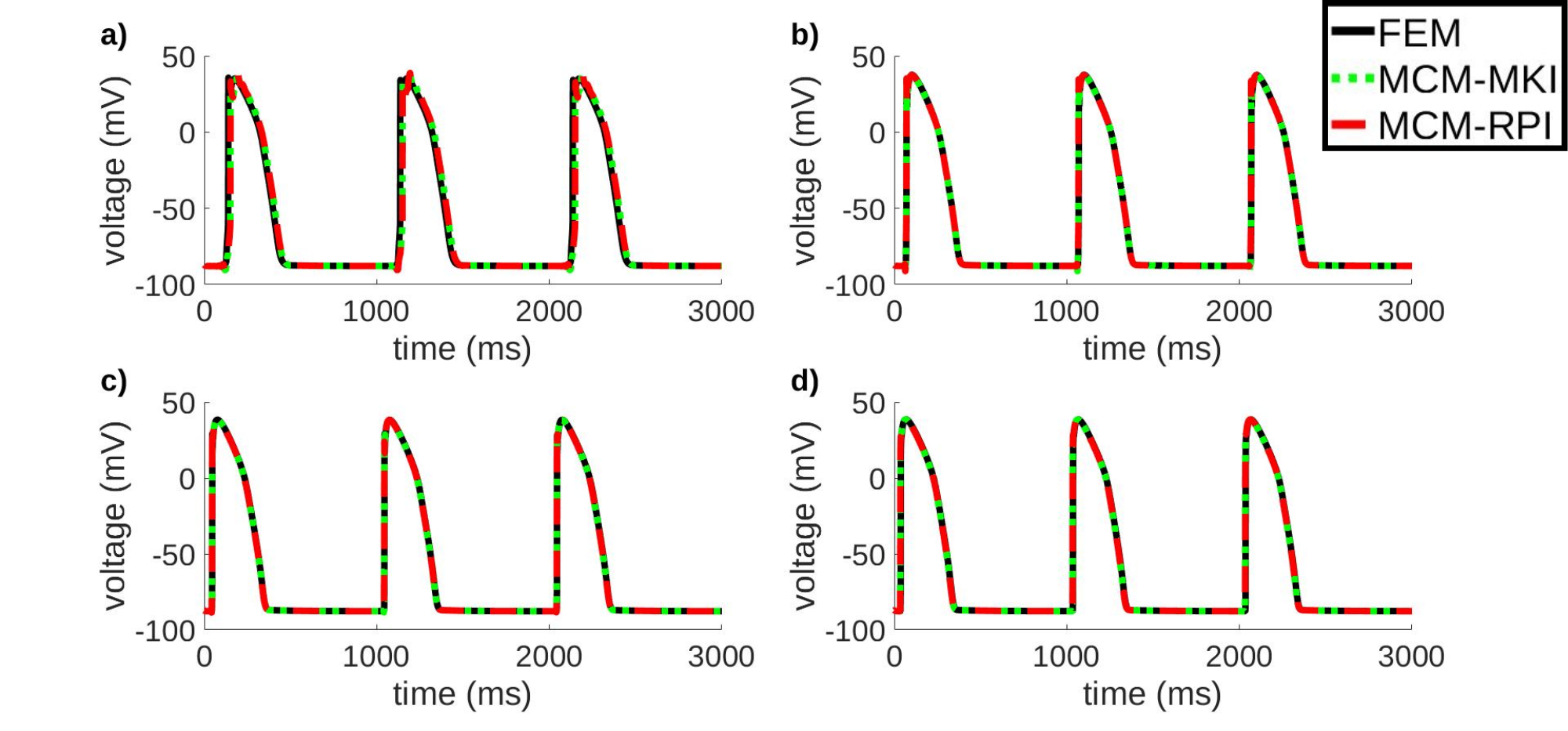}
    \caption{Voltage traces at the center ($x=2.5$ cm, $y=2.5$ cm) of a $5 \times 5$ cm ventricular tissue sheet calculated using FEM (continuous black), MCM with MKI (dotted green) and MCM with RPI (dashed red). The nodal discretization spacing is: a) $h=0.2$ cm, b) $h=0.1$ cm, c) $h=0.05$ cm, and d) $h=0.025$ cm.}
    \label{fig:2D_ap}
\end{figure}

We compare the MCM solution with RPI and MKI trial functions against FEM simulation results using bilinear isoparametric elements. We consider regular nodal discretizations and quadrilateral meshes with nodal spacing $h = \{0.2,\;0.1,\;0.05,\;0.025\}$ cm. The considered support domains in the meshfree approximation have size $s_d = \alpha_{sd} h$, with $\alpha_{sd} = 2.8$. The generated APs at the center of the tissue sheet ($x=2.5$ cm,\;$y=2.5$ cm) in the time interval $t = [0-3]$ s for the different nodal spacing values are shown in Figure \ref{fig:2D_ap}. We quantify the differences between MCM and FEM solutions in terms of mean transmembrane potential difference (TPD). Mean TPD between FEM and MCM with MKI trial functions was TPD = \{3.111,\;0.339,\;0.401,\;0.583\} mV while mean TPD between FEM and MCM with RPI trial functions was TPD = \{3.112,\;0.340,\;0.400,\;0.582\} mV for nodal spacing $h = \{0.2,\;0.1,\;0.05,\;0.025\}$ cm. The efficiency of each simulation is evaluated in terms of execution time in Figure \ref{fig:sim_time}\textcolor{medium_blue}{a}.

\subsection{Electrical propagation in a 3D tissue slab}

\begin{figure}[H]
    \centering
    \includegraphics[width=\textwidth]{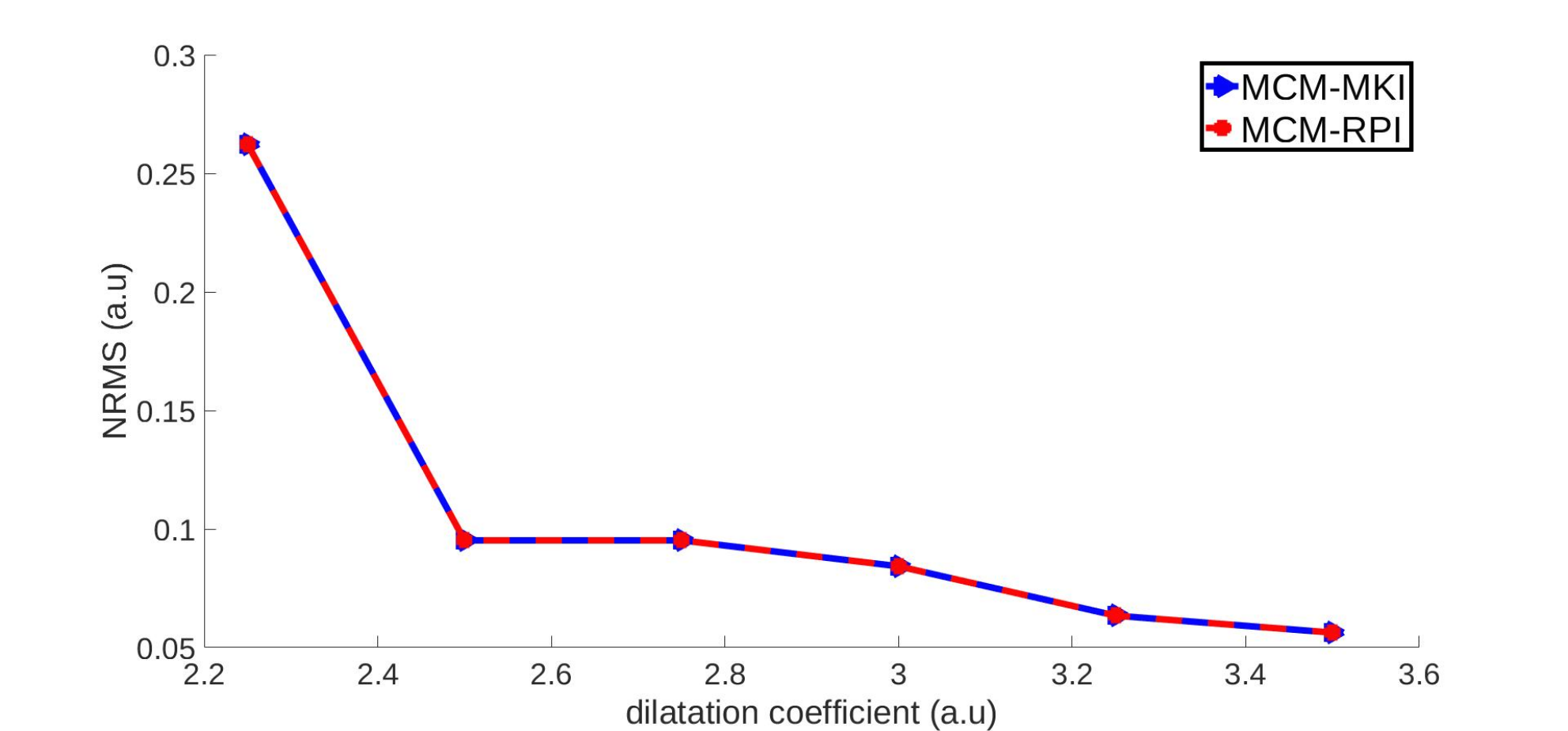}
    \caption{Convergence rate for varying support domain's dilatation coefficient $\alpha_{sd} \in \{2.25,\; 2.50,\; 2.75,\; 3.00,\; 3.25, \; 3.50\}$ for MCM with MKI (continuous blue) and RPI (slashed red) trial functions. Convergence is evaluated by computing the NRMS error with respect to a FEM simulation.}
    \label{fig:3D_conv}
\end{figure}

We investigate the effect of the support domain's dilatation coefficient $\alpha_{sd}$ by computing the normalized root mean square (NRMS) error between MCM and FEM solutions for a $3 \times 3 \times 3$ cm tissue slab. The tissue is assumed to be composed of epicardial ventricular cells. Stimuli of amplitude twice diastolic threshold, period $t_T = 1$ s and duration $t_d = 1$ ms are applied onto the left side of the tissue slab ($x=0$ cm). The tissue slab is discretized with $h = 0.05$ cm and varying dilatation coefficient $\alpha_{sd} \in= \{2.25,\; 2.50,\; 2.75,\; 3.00,\; 3.25, \; 3.50\}$. The NRMS error is computed using the formula:
\begin{equation} \label{nrms_error}
    \textrm{NRMS} = \frac{\sqrt{\displaystyle \frac{\displaystyle \sum_{I=1}^N \left( V_{MCM_I} - V_{FEM_I} \right)^2}{N}}}{\displaystyle \max_{I} V_{FEM_{I}}  - \displaystyle \min_{I} V_{FEM_{I}}}
\end{equation}

\noindent where $V_{MCM_I}$ and $V_{FEM_I}$ denote the transmembrane potential value at node $I$ computed with MCM and FEM, respectively. The NRMS error convergence plot is given in Figure \ref{fig:3D_conv}.

\begin{figure}[H]
    \centering
    \includegraphics[width=\textwidth]{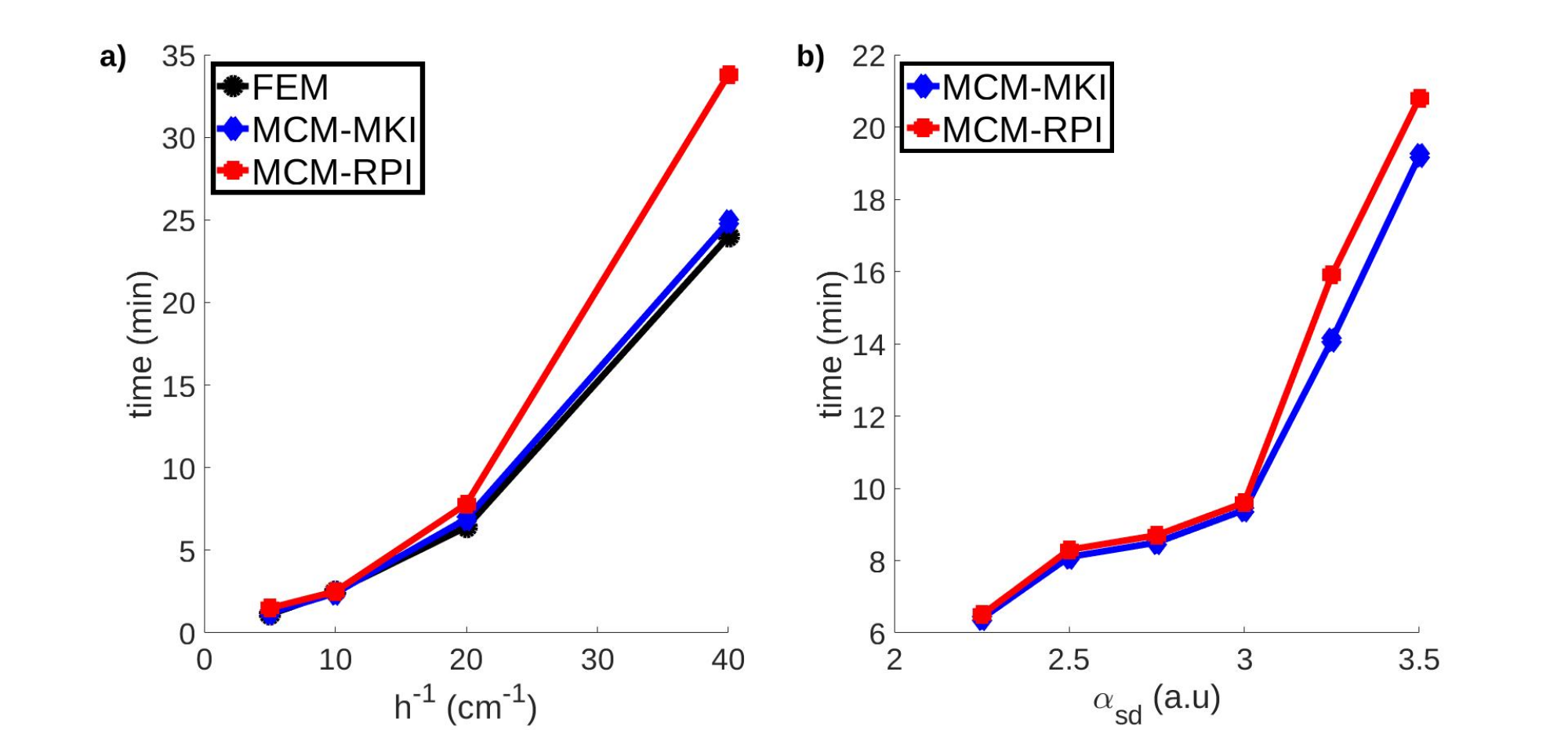}
    \caption{Execution time for a) varying spacing in a $5 \times 5$ cm tissue sheet and b) varying support dilatation coefficient in a $3 \times 3 \times 3$ cm tissue slab. Execution time is reported for FEM (black), MCM with MKI trial functions (blue) and MCM with RPI trial functions (red).}
    \label{fig:sim_time}
\end{figure}

The maximum NRMS error of the MCM solutions with RPI or MKI trial functions is obtained for $\alpha_{sd}=2.25$ and is equal to 0.262 and 0.263, respectively. The minimum NRMS error is obtained for $\alpha_{sd}=3.5$ and it is equal to 0.056 for RPI and 0.057 for MKI trial functions. The execution time for the simulations with varying dilatation coefficient is summarized in Figure \ref{fig:sim_time}\textcolor{medium_blue}{b}.

\subsection{Electrical propagation in a 3D biventricular geometry under left bundle branch block conditions} \label{subsec:lbba}

We simulate electrical propagation in a porcine cardiac biventricular model under healthy (baseline) and left bundle branch block (LBBB) conditions. The biventricular model is represented by a tetrahedral mesh (273,919 nodes, 1,334,218 elements). It is provided by the CRT-EPiggy19 challenge \cite{camara2019best,pop2020statistical} and is part of an LBBB dataset for experimental studies of cardiac resynchronization therapy \cite{rigol2013development,iglesias2016quantitative}. We compute the direction of the myocardial fibers using a rule-based method \cite{doste2019rule}. 

We define a diffusion coefficient value $d_0 = 0.002$ cm$^2$/ms in the longitudinal direction of the myocardial fibers and a transverse-to-longitudinal conductivity ratio $\rho = 0.2$. The electrophysiology of the ventricular myocardial tissue is represented using the O'Hara et al. model, as in previous examples. For a portion of connective tissue within the ventricles, the diffusion coefficient is reduced by a factor of 3 and the electrophysiology is modeled using the MacCannell et al. active fibroblast model \cite{maccannell2007mathematical}. 

To identify a set of points with the earliest activation across the ventricles,  we coupled the biventricular model with a network of Purkinje fibres generated using a fractal-tree generation algorithm \cite{costabal2016generating}. We applied stimuli with $t_d = 1$ ms, $t_T = 1$ s and amplitude twice the diastolic threshold onto the Purkinje-Myocardial Junctions (PMJs) identified from the terminal nodes of the Purkinje network. Electrical impulse propagation is simulated using MCM with RPI and MKI approximants as well as FEM.

The nodal support domains in MCM simulations are constructed using the nearest-neighbor approach with the 150 nearest nodes of each field node included in its support domain. We adopt this approach due to the irregular distribution of the nodes in the tetrahedral mesh. For such a distribution, constructing dilated support domains requires using a large value of the dilatation coefficient to avoid numerical instability, which leads to a very large number of neighboring nodes. Here, we use 150 nearest nodes to accurately capture the steep voltage gradients of the monodomain model.

\begin{figure}[H]
    \centering
    \includegraphics[width=\textwidth]{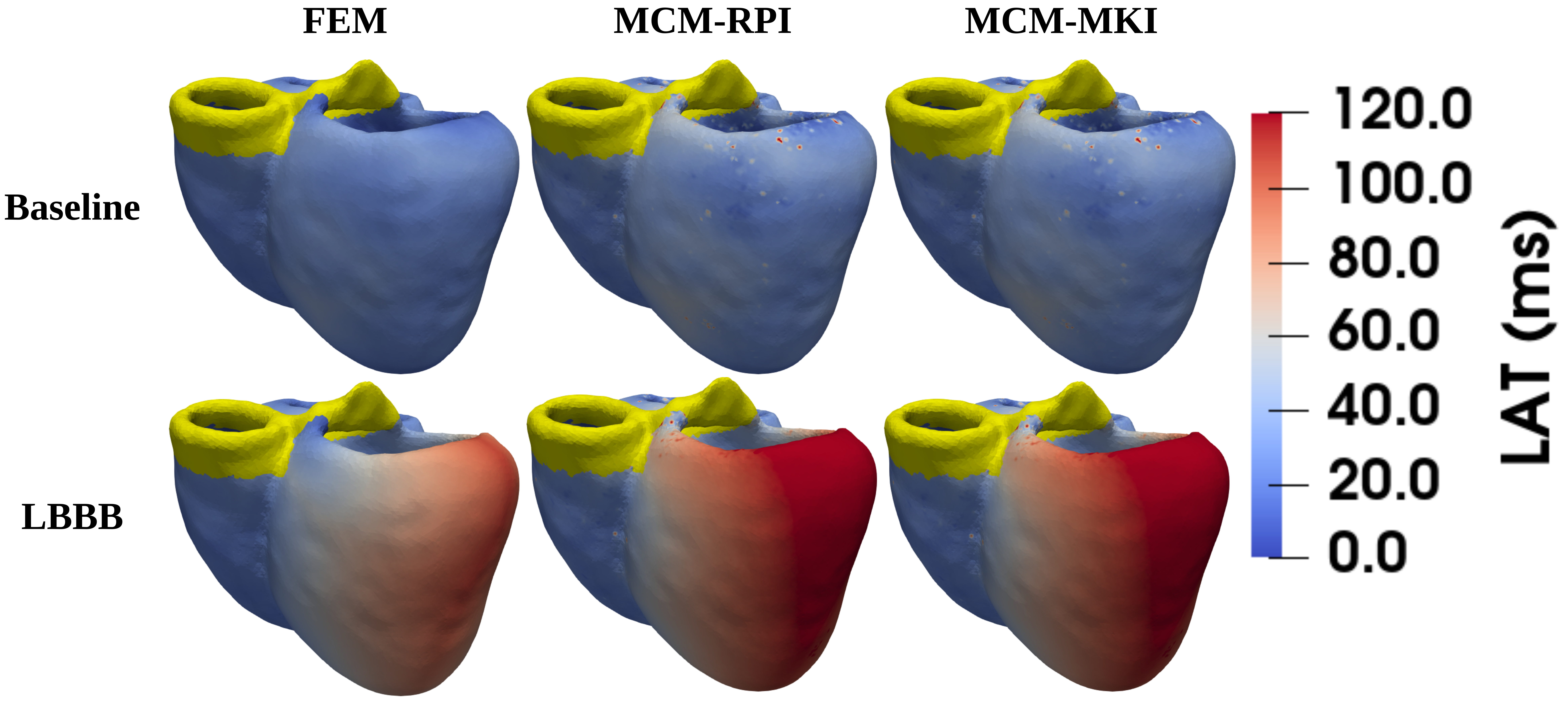}
    \caption{Local activation time maps for FEM (left), MCM-RPI (center), and MCM-MKI (right) simulations on a biventricular model. Connective tissue is represented in yellow.}
    \label{fig:lbba_lat}
\end{figure}

In Figure \ref{fig:lbba_lat}, we provide the MCM and FEM simulation results in terms of local activation time (LAT) maps for both baseline and LBBB conditions. 

The mean local activation time is computed for internal field nodes (LAT$_{in}$) and surface field nodes (LAT$_s$) separetely. Mean LAT$_{in}$ is found to be 23.2 ms for MCM-RPI and MCM-MKI under healthy conditions, while mean LAT$_{in}$ for FEM is 20.9 ms. Mean LAT$_{s}$ is found to be 27.3 ms for MCM-RPI and MCM-MKI and 21.4 ms for FEM. For LBBB conditions, mean LAT$_{in}$ is 47.9 ms for MCM-RPI and MCM-MKI and 40.4 ms for FEM. Mean LAT$_s$ is 51.9 ms for MCM-RPI and MCM-MKI and 41.1 ms for FEM. The execution time under healthy and LBBB conditions is 15.7 and 16.2 mins for MCM-RPI, 15.3 and 16.1 mins for MCM-MKI and 7.7 and 7.9 mins for FEM. For MCM-RPI and MCM-MKI methods, an additional time of 14.7 and 11.6 mins, respectively, is required for the calculation of the trial functions and their gradients. 

\subsection{Electrical propagation in a 3D biventricular geometry under myocardial infarction conditions}

We evaluate the MCM method for the simulation of myocardial infarction conditions using the biventricular anatomy of section \ref{subsec:lbba}. We introduce an infarct scar at the anterior wall of the left ventricle and we assume that the conductivity of the scarred tissue is zero. For these simulations, we do not model the border zone between the infarct and non-infarct areas. We compare the solution provided by MCM-RPI and MCM-MKI with the FEM solution in terms of LAT. The MCM is solved by representing the biventricular anatomy with a tetrahedral mesh, also used for FEM, as well as by using an immersed grid model \cite{mountris2020next}. 

\begin{figure}[H]
    \centering
    \includegraphics[width=\textwidth]{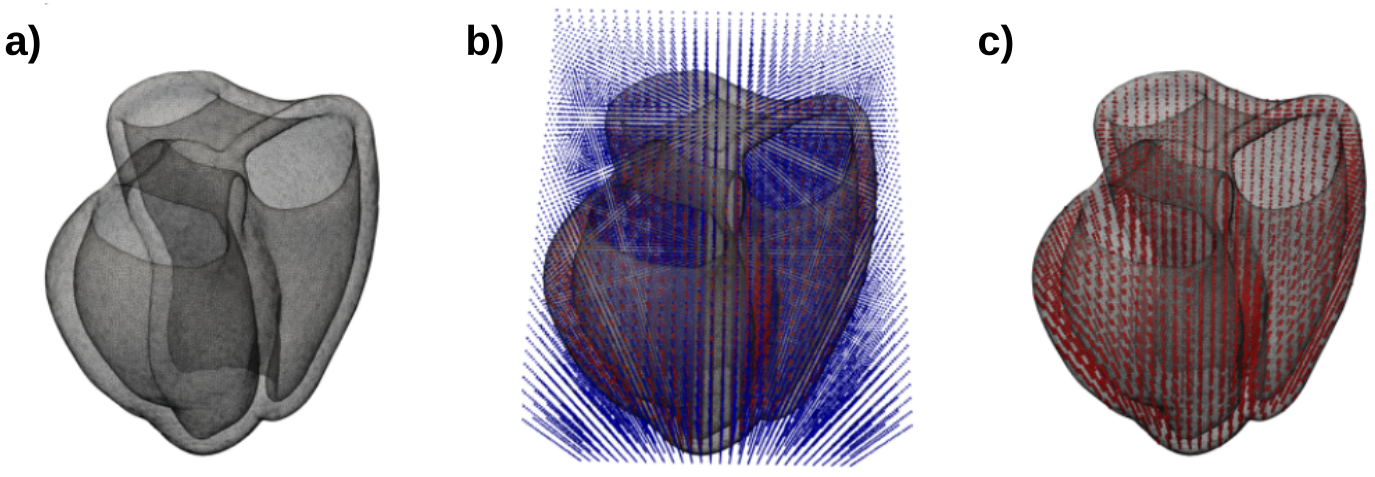}
    \caption{Immersed grid model generation. a) The boundary surface mesh is extracted from the tetrahedral biventricular mesh. b) Equidistant nodes with spacing twice the mean circumference of the triangles in the surface mesh are distributed inside its bounding box. c) Equidistant nodes located outside the surface mesh are discarded using a point containment test algorithm and the final immersed grid model is given by the union of the surface mesh nodes and the equidistant nodes enclosed by the surface.}
    \label{fig:immersed_grid}
\end{figure}

The immersed grid (238,519 nodes) is generated by distributing equidistant nodes inside the bounding box of the biventricular anatomy's boundary surface mesh. The nodal spacing of the immersed grid nodes is selected as twice the mean circumference of the triangles in the surface mesh, $h=0.9$ mm. A point containment test algorithm \cite{liu2010new} is used to discard nodes that span outside of the surface's boundary. Consequently, the final immersed grid is composed of the equidistant nodes inside the biventricular model and the nodes on the boundary surface mesh of the model (Figure \ref{fig:immersed_grid}).

\begin{figure}[H]
    \centering
    \includegraphics[width=\textwidth]{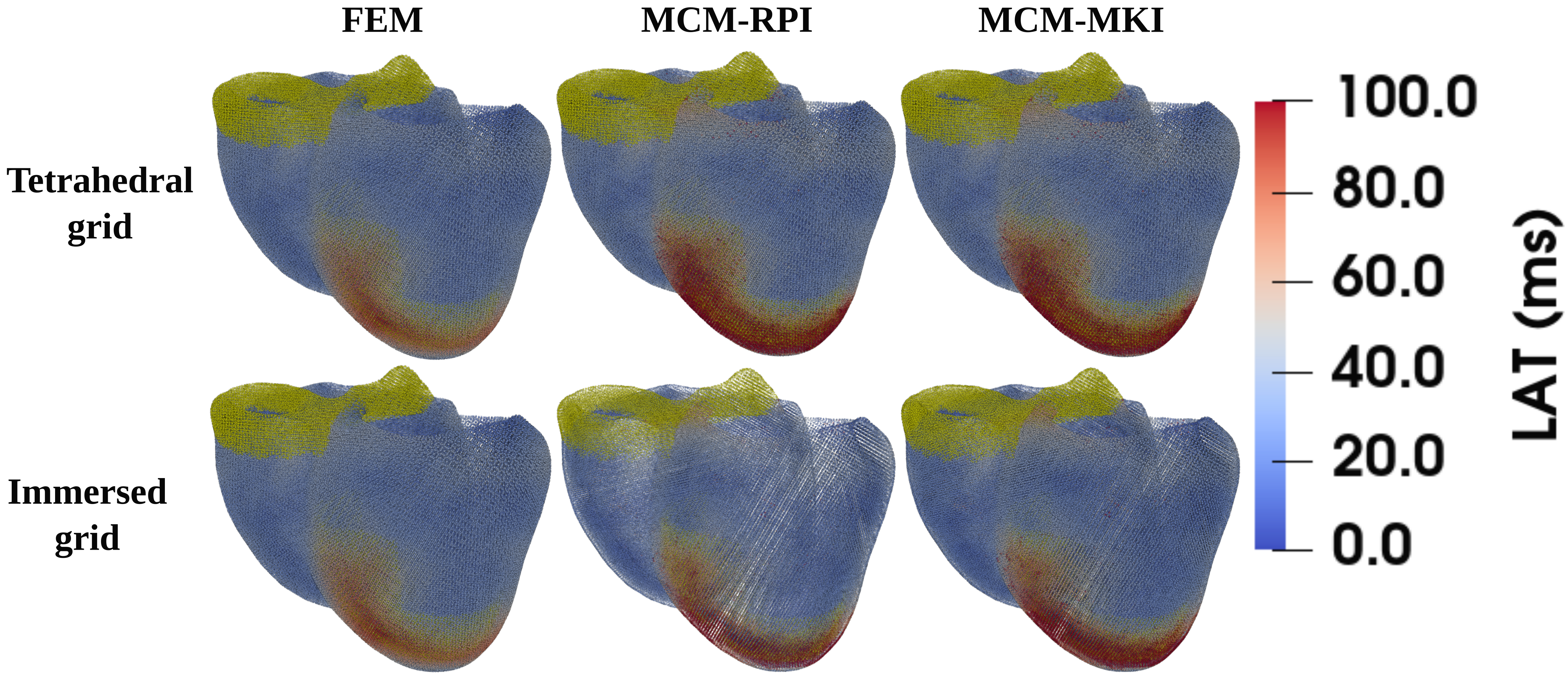}
    \caption{Local activation time maps for FEM (left), MCM-RPI (center), and MCM-MKI (right) simulations on an immersed grid of the biventricular model. Non-conductive connective tissue is represented in yellow.}
    \label{fig:scar_lat}
\end{figure}

We use the same tissue properties and electrical stimulation protocol as in section \ref{subsec:lbba} to generate the LAT maps for MCM-RPI, MCM-MKI and FEM for the tetrahedral and immersed grid models of the biventricular anatomy with scarred tissue. The LAT maps are shown in Figure \ref{fig:scar_lat}. For the tetrahedral model, the mean LAT$_{in}$ is 29.8 ms for MCM-RPI and MCM-MKI and 21.6 for FEM. Mean LAT$_s$ is 36.5 ms for MCM-RPI and MCM-MKI and 22.8 for FEM. The execution time is 16.1 mins for MCM-RPI, 15.3 mins for MCM-MKI and 7.8 mins for FEM. The required time to compute the RPI trial function and gradient is 14.5 mins and for the MKI trial function and gradient it is 11.6 mins. For the immersed grid model, the mean LAT$_{in}$ is 25.2 ms for MCM-RPI and MCM-MKI, while the mean LAT$_s$ 33.4 ms for MCM-RPI and MCM-MKI. The execution time is 13.8 mins for MCM-RPI and 13.5 mins for MCM-MKI. Additional time of 12.6 and 9.8 mins is required for the calculation of the RPI and MKI trial functions and gradients for the immersed grid model.

\section{Concluding remarks} \label{sec:conclusions}

In this study, we derived the mixed collocation method (MCM) to solve the monodomain model for the numerical simulation of cardiac electrophysiology. We considered two different interpolating trial functions, the radial point interpolation (RPI) and the moving Kriging interpolation (MKI). We solved several numerical examples in 2D and 3D domains comparing the MCM solution with a solution obtained by FEM. The accuracy of MCM solutions was very similar for both RPI and MKI approximants. However, MKI was found to be more efficient since, in contrast to RPI, it does not require polynomial enrichment. 

In all the numerical examples in 2D tissue sheets and 3D tissue slabs, good agreement was found between MCM and FEM solutions. A convergence analysis in 3D demonstrated that the MCM solution improves for larger support domains as the number of included collocation points is increased. For large scale problems, including the 150 nearest neighbours of each field node in its support domain was found to be an optimal choice balancing accuracy and memory footprint.

In large scale simulations of a biventricular anatomy model, MCM was able to produce LAT maps in good agreement with FEM under healthy and disease (LBBB, myocardial infarction) conditions. The largest differences in LAT between MCM and FEM were mainly found at the surface nodes of the model. This limitation was attributed to the negative effect of non-smooth changes of normal vectors in the imposition of Neumann boundary conditions. Introducing the immersed grid model approach, we were able to obtain improved results both in terms of accuracy and efficiency.

MCM is proved to be a promising alternative to FEM for cardiac electrophysiology simulation since its meshfree nature alleviates the need for the generation of a mesh and could thus allow fast model generation in a clinical setting.

\section*{Acknowledgements}
This work was supported by the European Research Council under grant agreement ERC-StG 638284, by Ministerio de Ciencia e Innovaci\'on (Spain) through project PID2019-105674RB-I00 and by European Social Fund (EU) and Arag\'on Government through BSICoS group (T39\_20R). Computations were performed by the ICTS NANBIOSIS (HPC Unit at University of Zaragoza).

\bibliographystyle{unsrt}  
\bibliography{main}

\end{document}